\documentclass[11pt]{article} 
\usepackage{latexsym} 
\usepackage{amsmath}
\usepackage{amssymb} 
\usepackage{amscd} 
\usepackage{graphicx}
\begin{document} 
\newtheorem{Th}{Theorem}[section]
\newtheorem{Cor}{Corollary}[section]
\newtheorem{Prop}{Proposition}[section]
\newtheorem{Lem}{Lemma}[section]
\newtheorem{Def}{Definition}[section]
\newtheorem{Rem}{Remark}[section]
\newtheorem{Ex}{Example}[section]
\newtheorem{stw}{Proposition}[section]


\newcommand{\bet}{\begin{Th}}
\newcommand{\ent}{\stepcounter{Cor}
   \stepcounter{Prop}\stepcounter{Lem}\stepcounter{Def}
   \stepcounter{Rem}\stepcounter{Ex}\end{Th}}


\newcommand{\bec}{\begin{Cor}}
\newcommand{\enc}{\stepcounter{Th}
   \stepcounter{Prop}\stepcounter{Lem}\stepcounter{Def}
   \stepcounter{Rem}\stepcounter{Ex}\end{Cor}}
\newcommand{\bep}{\begin{Prop}}
\newcommand{\enp}{\stepcounter{Th}
   \stepcounter{Cor}\stepcounter{Lem}\stepcounter{Def}
   \stepcounter{Rem}\stepcounter{Ex}\end{Prop}}
\newcommand{\bel}{\begin{Lem}}
\newcommand{\enl}{\stepcounter{Th}
   \stepcounter{Cor}\stepcounter{Prop}\stepcounter{Def}
   \stepcounter{Rem}\stepcounter{Ex}\end{Lem}}
\newcommand{\bef}{\begin{Def}}
\newcommand{\enf}{\stepcounter{Th}
   \stepcounter{Cor}\stepcounter{Prop}\stepcounter{Lem}
   \stepcounter{Rem}\stepcounter{Ex}\end{Def}}
\newcommand{\ber}{\begin{Rem}}
\newcommand{\enr}{
   \stepcounter{Th}\stepcounter{Cor}\stepcounter{Prop}
   \stepcounter{Lem}\stepcounter{Def}\stepcounter{Ex}\end{Rem}}
\newcommand{\bee}{\begin{Ex}}
\newcommand{\ene}{
   \stepcounter{Th}\stepcounter{Cor}\stepcounter{Prop}
   \stepcounter{Lem}\stepcounter{Def}\stepcounter{Rem}\end{Ex}}
\newcommand{\Proof}{\noindent{\it Proof\,}:\ }

\newcommand{\EE}{\mathbf{E}}
\newcommand{\QQ}{\mathbf{Q}}
\newcommand{\R}{\mathbf{R}}
\newcommand{\C}{\mathbf{C}}
\newcommand{\ZZ}{\mathbf{Z}}
\newcommand{\NN}{\mathbf{N}}
\newcommand{\PP}{\mathbf{P}}
\newcommand{\uuu}{{u}}
\newcommand{\xxx}{{x}}
\newcommand{\aaa}{{a}}
\newcommand{\bbb}{{b}}
\newcommand{\AAA}{\mathbf{A}}
\newcommand{\BBB}{\mathbf{B}}
\newcommand{\ccc}{{c}}
\newcommand{\iii}{{i}}
\newcommand{\jjj}{{j}}
\newcommand{\kkk}{{k}}
\newcommand{\rrr}{{r}}
\newcommand{\FFF}{{F}}
\newcommand{\yyy}{{y}}
\newcommand{\ppp}{{p}}
\newcommand{\qqq}{{q}}
\newcommand{\nnn}{{n}}
\newcommand{\vvv}{{v}}
\newcommand{\eee}{{e}}
\newcommand{\fff}{{f}}
\newcommand{\www}{{w}}
\newcommand{\0}{{0}}
\newcommand{\lon}{\longrightarrow}
\newcommand{\ga}{\gamma}
\newcommand{\pa}{\partial}
\newcommand{\QED}{\hfill $\Box$}
\newcommand{\id}{{\mbox {\rm id}}}
\newcommand{\Ker}{{\mbox {\rm Ker}}}
\newcommand{\Tan}{{\mbox {\rm Tan}}}
\newcommand{\grad}{{\mbox {\rm grad}}}
\newcommand{\ind}{{\mbox {\rm ind}}}
\newcommand{\rot}{{\mbox {\rm rot}}}
\newcommand{\diver}{{\mbox {\rm div}}}
\newcommand{\Gr}{{\mbox {\rm Gr}}}
\newcommand{\rank}{{\mbox {\rm rank}}}
\newcommand{\ord}{{\mbox {\rm ord}}}
\newcommand{\sign}{{\mbox {\rm sign}}}
\newcommand{\Spin}{{\mbox {\rm Spin}}}
\newcommand{\Symp}{{\mbox {\rm Sp}}}
\newcommand{\Int}{{\mbox {\rm Int}}}
\newcommand{\codim}{{\mbox {\rm codim}}}
\def\mod{{\mbox {\rm mod}}}
\newcommand{\qed}{\hfill $\Box$ \par}

\newcommand{\dint}[2]{{\displaystyle\int}_{{\hspace{-1.9truemm}}{#1}}^{#2}}


\title{Singularities of parallels to tangent developable surfaces} 

\author{Goo \textsc{Ishikawa}\thanks{Department of Mathematics, Faculty of Sciences, Hokkaido University, Sapporo 060-0810, Japan. 
e-mail: 
ishikawa@math.sci.hokudai.ac.jp}
}

\renewcommand{\thefootnote}{\fnsymbol{footnote}}
\footnotetext{The author is supported by KAKENHI no. 19K03458. }

\date{ }

\maketitle

\begin{abstract}
It is known that the class of developable surfaces which have zero Gaussian curvature in three dimensional Euclidean space 
is preserved by the parallel transformations.  
A tangent developable surface is defined as a ruled developable surface by tangent lines to a space curve and 
it has singularities at least along the space curve, called the directrix or the the edge of regression. 
Also the class of tangent developable surfaces are invariant under the parallel deformations. 
In this paper the notions of tangent developable surfaces and their parallels are naturally generalized 
for frontal curves in general in Euclidean spaces of arbitrary dimensions.  
We study singularities appearing on parallels to tangent developable surfaces of frontal curves and 
give the classification of generic singularities on them for frontal curves in 3 or 4 dimensional Euclidean spaces. 
\end{abstract}

\section{Introduction}

Given a surface in $\R^3$, 
a \lq\lq parallel surface" or a "parallel" to the surface is simply defined as a surface which have the common family of normal affine lines with the original surface. In fact, given a surface $(u, v) \mapsto f(u, v)$ with a unit normal $\nu(u, v)$, 
its parallels are given by the surfaces $f(u, v) + r \nu(u, v)$ with the parameter $r \in \R$. 

In general a line congruence, i.e., a two-dimensional family of affine lines in $\R^3$ is called a {\it system of rays} if it forms a (possibly singular) Lagrangian surface in the space of affine lines in $\R^3$ \cite{Arnold1, Arnold}. The condition is equivalent to that the family is expressed in a parametric form as 
$$
(u, v) \mapsto f(u, v) + r \,\nu(u, v) \ (r \in \R)
$$
for some \lq\lq frontal" surface $(u, v) \mapsto f(u, v) \in \R^3$ and 
its unit normal field $(u, v) \mapsto \nu(u, v) \in S^2$ so that 
$(u, v) \mapsto (f(u, v), \nu(u, v)) \in \R^3\times S^2$ is 
Legendrian, i.e. it is an integral map for the standard contact structure on 
the unit tangent bundle $\R^3\times S^2$ of the Euclidean space $\R^3$. 

Parallel deformations are regarded as an important and interesting transformations of surfaces 
in differential geometry.  
For instance, constant mean curvature surfaces 
and positive constant Gaussian curvature surfaces are related as parallels to each other. 
Parallels are studied deeply from geometric point of view of differential geometric surface theory  (see \cite{Spivak} p.185, \cite{IL} p.225 for instance). 
It is characteristic also that, in the process of taking parallels, the surfaces may have 
singularities with a geometric origin. Singularites which appear in parallels are studied 
deeply from singularity theory viewpoint \cite{Bruce, FH, IFRT}. 

A surface in Euclidean 3-space is called {\it developable} if it is isometric to the plane and the class of developable surfaces. 
The class is characterized also as the surfaces with zero Gaussian curvature.  
Developable surfaces are roughly classified into cylinders, cones and tangent developable surfaces. 
A tangent developable surface is defined as a ruled surface by tangent lines to a space curve and 
has singularities at least along the space curve, called the directrix or the edge of regression (\cite{Ishikawa99}). 
Note that the class of developable surfaces is preserved under the parallel transformation. 
Moreover any conical surface can be regarded as a tangent surface of the \lq\lq vertex" i.e. a constant curve endowed with 
a tangential framing along the constant curve, while the class of cylindrical surfaces are preserved under the parallel transformations. 

Then we are led naturally to ask which singularities appear on parallels to 
tangent developable surfaces of curves and the possibility to classify such generic singularities. 

The notion of tangent developable surfaces is generalized to frontal curves. 
A possibly singular curve $t \mapsto f(t)$ is called a {\it frontal} curve if there exists a family of unit vector field $\tau(t)$ 
such that $f'(t)$ is a scalar multiple of $\tau(t)$ for any $t$. 
A possibly singular surface $(u, v) \mapsto F(u, v)$ is called a {\it frontal} surface if 
there exists a family of unit vector field $\nu(u, v)$ such that $\nu(u, v)$ is orthogonal to both 
$F_u(u, v)$ and $F_v(u, v)$ for any $(u, v)$ (see \S \ref{Parallels to normally flat frontals}). 

Under the context mentioned above, we are going to study in this paper 
parallel frontals for the case when a system of rays in $\R^3$ is \lq\lq degenerate", i.e. when 
the corresponding Lagrangian surface in the space of affine lines in $\R^3$ is expressed as 
$$
(t, s) \mapsto f(t) + s \tau(t) + r \nu(t) \ (r \in \R), 
$$
where $f(t)$ is a frontal curve in $\R^3$, $\tau(t)$ is a unit tangential frame along $f$ and 
$\nu(t)$ is a unit normal to the surface $(t, s) \mapsto f(t) + s \tau(t)$ of $f$, 
the tangent developable surface, or simply, the tangent surface of $f$. 
Here the degeneracy means that the direction vector field $\nu(t)$ for the normal lines is independent of $s$. 

The same idea is applied to frontal hypersurfaces or systems of rays in $\R^{n+1}$ and also 
even to $n$-dimensional frontal submanifolds or systems of \lq\lq multi-dimensional rays" 
in $\R^{n+p}$ when the frontal submanifolds are assumed to be \lq\lq normally flat"
(see \S \ref{Parallels to normally flat frontals} and \cite{Ishikawa20b}).

Then, in fact, we clarify that any parallels to the tangent surface of a curve turns to be again a tangent surface of some curve, 
under a mild condition, and 
thus the class of tangent developable surfaces of curves is in fact preserved by the surface transformation 
by taking parallels. Furthermore we show that the observation remains true for frontal 
curves not only in $\R^3$ but also in $\R^{1+p}, p \geq 2$ in general 
(see \S \ref{Parallels to normally flat frontals} and \cite{Ishikawa20b}). 
Actually we obtain the explicit formula for directrixes which produce parallels by taking tangent lines 
in terms of differential geometric data on the original curve (Theorem \ref{parallel-of-tangent}). 
In this paper we proceed to classify singularities which appear in parallel surfaces to the 
tangent surfaces of generic frontal curves in $\R^3$ and $\R^4$ (Theorem \ref{classification-p=2} and 
Theorem \ref{classification-p=3}).

In the next section \S \ref{Parallels to normally flat frontals} we make clear the framework of this paper, introducing 
Theorem \ref{parallel-of-tangent}, and 
in \S \ref{Classification of generic singularities} we formulate our classification results in this paper (Theorems \ref{classification-curve-type}, 
\ref{classification-p=2} and \ref{classification-p=3}). 
In \S \ref{Parallels to tangent surfaces of frontal curves} we prove Theorem \ref{parallel-of-tangent} briefly 
and study singularities of parallels to tangent surfaces to 
curves in $\R^{n+p}$ in general to show Theorems \ref{classification-curve-type} 
\ref{classification-p=2} and \ref{classification-p=3}. 
In \S \ref{The cuspidal swallowtail surfaces} and in \S \ref{The swallowtail and its openings} 
we study the classification problem of some singularities which appear in the list of Theorem \ref{classification-p=3}. 

\

All manifolds and mappings in this paper are assumed of class $C^\infty$ unless otherwise stated.

\section{Parallels to normally flat frontals} 
\label{Parallels to normally flat frontals} 

Let $M$ be an $n$-dimensional manifold and $\R^{n+p}$ the $(n+p)$-dimensional Euclidean space. 

A map-germ $f : (M, a) \to \R^{n+p}$ at a point $a \in M$ is 
called a {\it frontal-germ} 
if there exists an orthonormal frame $\nu_1, \nu_2, \dots, \nu_p : (M, a) \to T\R^{n+p}$ along $f$ 
such that 
$$
\nu_i(u) \cdot f_*(T_u M) = 0
$$ 
for any $u \in M$ nearby $a$ and for any $i = 1, 2, \dots p$ (\cite{Ishikawa18}). 
Hear $\cdot$ means Euclidean inner product and $f_* : TM \to T\R^{n+p}$ the differential of $f$. 
Then we consider the subbundle $N_f$ generated by $\nu_1, \nu_2, \dots, \nu_p$ as above 
in the pull-back bundle $f^*(T\R^{n+p})$. We call $N_f$ the {\it normal bundle} 
of the frontal-germ $f$. 

A mapping $f$ is called a {\it frontal} if the germ of $f$ at any point $a \in M$ is a frontal-germ. 

We remark that, in general, $N_f$ is not uniquely determined if $f$ is very degenerate. 
We call a mapping $f$ defined on $M$ {\it proper} if 
the singular locus $\Sigma(f)$, i.e. the non-immersive locus of $f$, is nowhere dense in $M$. 

For a proper frontal $f : M \to \R^{n+p}$, 
the normal bundle $N_f$ is uniquely determined on $M$, because it is determined as the orthogonal complement to 
$f_*(TM)$ which is regarded as a subbundle of $f^*(T\R^{n+p})$ over the regular locus $M \setminus \Sigma(f)$ of $f$ (\cite{Ishikawa18} Proposition 6.2). 
Then, in turn, the {\it tangent bundle} of the frontal $f$ is defined as the orthogonal complement to $N_f$ in 
$f^*(T\R^{n+p})$ over the whole $M$. 
Thus the pull-back bundle $f^*(T\R^{n+p})$ is decomposed into the 
sum $T_f \oplus N_f$ by the tangent bundle $T_f$ of rank $n$ 
and the normal bundle $N_f$ of rank $p$. 

A frontal $f$ is called {\it normally flat} (resp. {\it tangentially flat})  
if the induced connection on $N_f$ (resp. the induced connection on $T_f$) 
from the Euclidean connection on $f^*(T\R^{n+p})$ is flat (\cite{KN, Terng}). 
If $f$ is normally flat, then there exists an orthonormal frame 
$\{ \nu_1, \dots, \nu_p\}$ along $f$ such that the covariant derivative of each $\nu_i$ 
by any vector field over $M$ belongs to $T_f$. 
We call such a frame a {\it normally parallel orthonormal frame} or {\it Bishop frame} of the normally flat frontal (see \cite{Bishop}). 
Normally parallel sections of $N_f$ form a $p$-dimensional vector space and 
any normally parallel section $\nu$ of $N_f$ is written as a linear combination $\nu = \sum_{i=1}^p r_i\nu_i, 
(r_1, \dots, r_p \in \R)$ (See \cite{Ishikawa20b} Lemma 3.1). A normally parallel vector field $\nu$ 
is said to be {\it generic} if $r_1, \dots, r_p$ are taken to be generic. 

Frontal curves and hypersurfaces are normally flat. 
In \cite{Ishikawa20b}, we define, 
for a given normally flat frontal in a Euclidean space, 
its parallel frontals using normally flat frame (Bishop frame). 
Moreover it is shown that the tangent surface of a frontal curve which is generated by tangent lines 
to the curve is normally flat, provided it is a frontal, 
and every parallel to the tangent surface turns to be right equivalent to 
the tangent surface of a frontal curve, called the {\it directrix} or the {\it edge of regression}. 

Let $f : M \to \R^{n+p}$ be a normally flat frontal, and 
$\{ \nu_1, \dots, \nu_p\}$ a Bishop frame of $N_f$ defined over $M$. Then we define 
the {\it parallels} to $f$ by 
$$
{\textstyle 
{\mbox{\rm P}}_\nu(f)(t) = {\mbox{\rm P}}_r(f)(t) := f(t) + \nu(t), 
{\mbox{\rm where }} \nu(t) = \sum_{i=1}^p r_i \nu_i(t), \   (r \in \R^p). 
} 
$$
Then the parallels are normally flat and have the same Bishop frame with $f$. 

Parallels are basic and interesting objects to be studied in the cases of both hypersurfaces and curves 
\cite{FH, HT}. 

Let $f : M \to \R^{n+p}$ be a frontal. 
If $\{ \tau_1, \dots, \tau_n\}$ is a frame of the tangent bundle $T_f$ over $M$, then 
the {\it tangent map} ${\mbox{\rm Tan}}(f) : M\times \R^n \to 
\R^{n+p}$ is defined by
$$
{\textstyle 
{\mbox{\rm Tan}}(f)(u, s) = f(u) + \sum_{i=1}^n s_i \tau_i(u), \quad u \in M, s \in \R^n. 
}
$$
The right equivalence class of ${\mbox{\rm Tan}}(f)$ is independent of the choice of 
$\{ \tau_1, \dots, \tau_n \}$. 
When $n = 1$, then the tangent map is called a tangent surface. 

Then we have that tangent surfaces of frontal curves are normally and tangentially flat: 

\bet
\label{tangent-nomally-flat}
{\rm (\cite{Ishikawa20b}, Theorem 2.13)} \ 
Let $I$ be an interval and 
$f : I \to \R^{1+p}$ a frontal curve.  
Suppose the tangent surface $\Tan(f) : I\times\R \to \R^{1+p}$ is a proper frontal, i.e. a frontal with nowhere dense singular points. 
Then $\Tan(f)$ is a normally and tangentially flat frontal. 
\ent

Let $f : I \to \R^{1+p}$ be a frontal curve. A point $t \in I$ is called an {\it inflection point} 
of $f$ if there exists a unit frame $\tau : I \to S^p (\subset \R^{1+p})$ of $T_f$ such that 
$t$ is a singular point of $\tau$, regarded as a curve in the sphere $S^p$. 

Suppose $F = \Tan(f) : I\times \R \to \R^{1+p}$ is a proper frontal and 
$T_F = T_{\Tan(f)}$ is the tangent bundle of $\Tan(f)$. 
Let $\tau$ be a unit section of $T_f$ regarded as a unit section of $T_{\Tan(f)}$ along $I \times\{ 0\}$. 
Note that $\Tan(f)(t, 0) = f(t)$. 
Take an orthonormal frame $\{ \tau, \mu\}$ of $T_{\Tan(f)}$ along $I\times\{ 0\}$. 
Then there exists unique function $\kappa : I \to \R$ with $\tau'(t) = \kappa(t)\mu(t)$, 
which is called the {\it curvature function} of the frontal curve $f$ with respect 
the frame $\{ \tau, \mu\}$. Remark that $t_0 \in I$ is an inflection point of $f$ if 
and only if $\kappa(t_0) = 0$. 
Take a normally parallel orthonormal frame $\{ \nu_1, \dots, \nu_{p-1}\}$ of 
$N_{\Tan(f)}$ along $I\times\{ 0\}$. Then there exist functions 
$\ell_1, \dots, \ell_{p-1} : I \to \R$ uniquely, which is called {\it torsion functions} of $f$, 
such that $\mu'(t) = - \kappa(t)\tau(t) + 
\sum_{i=1}^{p-1} \ell_i(t)\nu_i(t)$. 
The frame $\{ \nu_1, \dots, \nu_{p-1}\}$ extends to a 
normally parallel orthonormal frame of $N_{\Tan(f)}$ over $I\times\R$ 
as $\nu_i(t, s) = \nu_i(t)$, by the constancy of $N_{\Tan(f)}$ along $\{ t\}\times \R$. 

Now let us consider parallels ${\mbox{\rm P}}_r(\Tan(f))$ of $\Tan(f)$ 
generated by $\{ \nu_1, \dots, \nu_{p-1}\}$ of $N_{\Tan(f)}$. 
Then we have that parallels to tangent surfaces of frontal curves are again tangent surfaces of frontal curves as follows: 

\bet
\label{parallel-of-tangent}
{\rm  (\cite{Ishikawa20b}, Theorem 2.17)} \ 
Let $p \geq 2$. 
Let $f : I \to \R^{1+p}$ be a frontal curve without inflection point, 
and suppose $\Tan(f) : I\times \R \to \R^{1+p}$ is a proper frontal. 
Then any parallel ${\mbox{\rm P}}_r(\Tan(f))$ 
to $\Tan(f)$ is right equivalent to the tangent surface $\Tan(g)$ 
for a frontal curve $g : I \to \R^{1+p}$. 
In fact $g$ is given by 
$$
g(t) = f(t) + \sum_{i=1}^{p-1} r_i\left(\frac{\ell_i(t)}{\kappa(t)}\,\tau(t) + \nu_i(t)\right). 
$$
If $f'(t) = a(t)\tau(t)$, then the velocity vector of $g$ is given by 
$$
g'(t) = 
\left( a(t) + \sum_{i=1}^{p-1} r_i\left\{ \frac{\ell_i(t)}{\kappa(t)}\right\}'\right)\tau(t). 
$$
\ent 

\bee
{\rm 
Let $f : (\R, 0) \to \R^3$ be a curve defined by 
$f(t) = (t, \frac{t^3}{6}, \frac{t^4}{24})$. 
Then $\tau(t) = (1, \frac{t^2}{2}, \frac{t^3}{6})$ and $f$ has an inflection point at $t_0 = 0$. 
The tangent surface is parametrized as  
$$
{\textstyle 
F(t, s) = (t+s, \ \frac{t^3}{6} + s\frac{t^2}{2}, \ \frac{t^4}{24} + s\frac{t^3}{6}).  
}
$$
The unit normal to $F$ is given by 
$$
{\textstyle 
\nu(t, s) = (\nu_1(t), \nu_2(t), \nu_3(t)) = \frac{12}{\sqrt{t^6+36t^2+144}}(\frac{t^3}{12}, -\frac{t}{2}, 1). 
}
$$
Then the parallels $F_r(t, s) = F(t, s) + r\nu(t)$ are never right equivalent to $\Tan(f_r)$ 
for any family of frontal curves $f_r : I \to \R^3$.

The following bifurcation of fronts is realized by the parallel deformation of a tangent developable surface (Mond surface) of a curve of type $(1, 3, 4)$. 
\vspace{-0.5truecm}

 \begin{center}  
   \includegraphics[width=10truecm, height=3.5truecm, clip, 
   ]{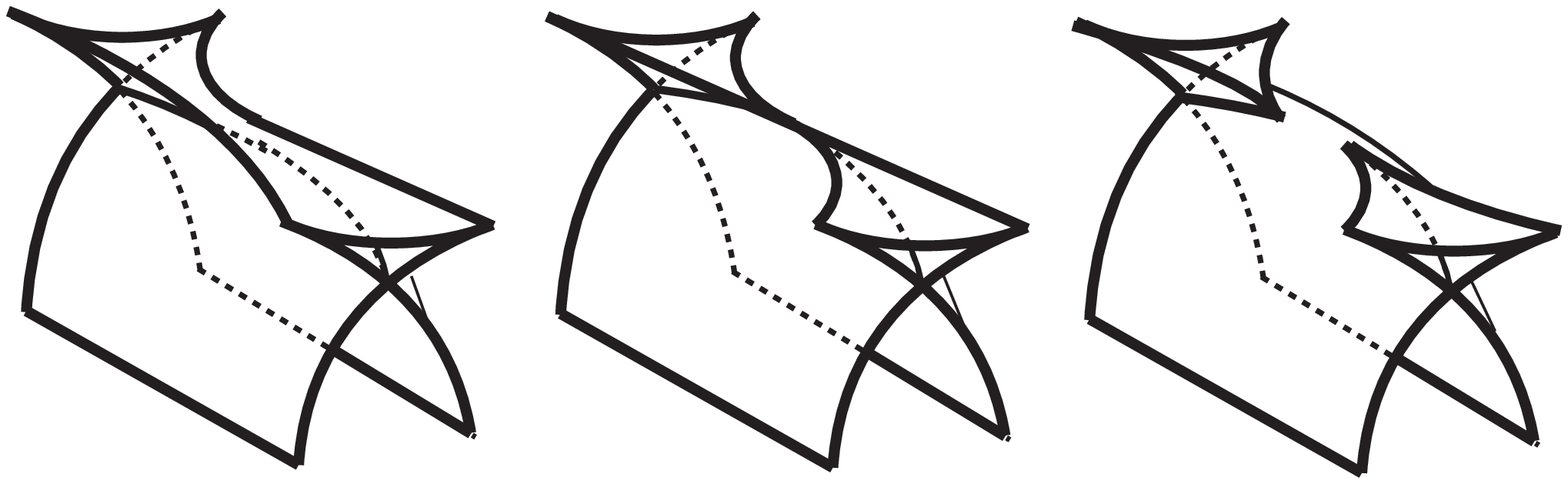}
   \\
   Parallel deformation of a Mond surface
\end{center}

}
\ene

\bee
{\rm 
Conical singularities can be obtained as tangent developable surfaces. 
Let $\tau : I \to S^2 \subset \R^3$ be a spherical curve such that $\kappa(t)$ does not take $0$ and 
$\ell_1(t)/\kappa(t)$ is not a constant. Set $a = - (\ell_1(t)/\kappa(t))'$ and 
take a frontal curve $f : I \to \R^3$ satisfying $f'(t) = a(t)\tau(t)$. Then the parallel 
$P_1(\Tan(f)) : I\times \R \to \R^3$ to $\Tan(f)$ degenerates to a conical surface. 
}
\ene

\section{Classification of generic singularities}
\label{Classification of generic singularities}

To give the classification results in this paper for singularities of parallels to 
tangent surfaces of frontal curves, we recall several necessary notions. 

We define, for a curve $f : I \to \R^{1+p}$, the $(n+p)\times k$-Wronskian matrix of $f$ by 
$$
W_k(f)(t) := 
\left( \dfrac{df}{dt}(t), \dfrac{d^2f}{dt^2}(t), \dfrac{d^3f}{dt^3}(t), \cdots, \dfrac{d^kf}{dt^k}(t)\right), 
$$
for any positive integer $k$. Note that the rank of $W_k(f)(t)$ is independent of the choices of 
local coordinate on $I$ and affine coordinates of $\R^{n+p}$. 

\bef
{\rm 
Let $a_1, a_2, \dots, a_{1+p}$ be an strictly increasing sequence of positive integers, 
$1 \leq a_1 < a_2 < \cdots < a_{1+p}$. 
A curve $f : I \to \R^{1+p}$ is said to be of type $(a_1, a_2, \dots, a_{1+p})$ at $t \in I$ if 
$a_i = \min\{ k \mid \rank(W_k(f)(t)) = i\}, i = 1, 2, \dots, 1+p$. 
}
\enf

By Theorem \ref{parallel-of-tangent}, the singularity of 
a parallel of $\Tan(f)$ are described by the tangent surface $\Tan(g)$ of a curve $g$, provided $f$ has no
inflection points. 

\bef
{\rm
We call a frontal curve $f : I \to \R^{1+p}$ {\it generic} if there exists an open dense subset in the space 
of maps $I \to \R\times S^p$ with respect to Whitney $C^\infty$ topology such that 
$f(t) = a(t)\tau(t)$ for some $(a, \tau)$ belonging to the open dense subset. 
}
\enf

For the basic notions of singularity theory of mappings, see \cite{GG, AGV} for instance. 

\bep
\label{generic-frontal-curves} 
Generic frontal curves $I \to \R^{1+p}$ have no inflection points, provided $p \geq 2$. 
In fact a generic frontal curve is, at each point of $I$, of type $(1, 2, \dots, p, 1+p), 
(1, 2, \dots, p, 2+p)$ or $(2, 3, \dots, 1+p, 2+p)$. 
\enp

Note that generic plane frontal curves $I \to \R^2$, $p = 1$, may have inflection points, 
which are of type $(1, 3)$. 

\

Now we give the classification of types of singularities appearing 
in directrixes of parallels to tangent surfaces of generic frontal curves 
generated by its generic normally parallel normal fields. 

\bef
{\rm 
A frontal curve $f : I \to \R^{1+p}$ is called {\it generic} if $f'(t) = a(t)\tau(t)$ for a mapping 
$(a, \tau) : I \to \R\times S^p$ which belongs to an open dense subset in the space $C^\infty(I, \R\times S^p)$ 
of $C^\infty$ mappings with respect to $C^\infty$-topology (\cite{GG}). 
}
\enf

\bet
\label{classification-curve-type} 
Let $p \geq 2$. 
Let $f : I \to \R^{1+p}$ be a generic frontal curve and 
$\Tan(f) : I \times \R \to  \R^{1+p}$ the tangent surface of $f$. 
Let $\nu : I\times \R \to \R^{1+p}$ be a generic normally parallel normal field along $\Tan(f)$. 
Then the type of the directrix $g : I \to \R^{1+p}$ 
of the parallel $\Tan(f) + r\nu : I\times \R \to \R^{1+p}$ for any $(t, r) \in I\times\R$ 
is given by 
$$
\begin{array}{|c|c|}
\hline
{\mbox {\rm type}} & {\mbox{\rm codimension}} 
\\
\hline
(1, 2, \dots, p, 1+p) & 0
\\
(1, 2, \dots, p, 2+p) & 1
\\
(2, 3, \dots, 1+p, 2+p) & 1
\\
(2, 3, \dots, 1+p, 3+p) & 2
\\
(3, 4, \dots, 2+p, 3+p) & 2
\\
\hline
\end{array}
$$
Each number in the second column means the codimension of the locus of the type on the $(t, r)$-plane $I\times\R$. 
\ent

It is known that, in several cases, singularities of tangent surfaces $\Tan(g)$ are described by types of $g$ 
(\cite{Ishikawa99, Ishikawa12}). 

We give in this paper the exact generic classification of parallels to tangent developables of frontal curves 
in $\R^3$ and in $\R^4$. 

\

Let $p = 2$. Then the possible types are 
$(1, 2, 3), (1, 2, 4), (2, 3, 4), (2, 3, 5)$ and $(3, 4, 5)$ by Theorem \ref{classification-curve-type}. 
Therefore we have

\bet
\label{classification-p=2} 
{\rm (Classification of singularities of generic parallels to tangent maps of frontal curves in $\R^3$.)} \
Let $p = 2$. 
Any singular point of parallels $f(t) + s\tau(t) + r\nu(t)$ along a unit normal $\nu$ 
to the tangent surface $f(t) + s\tau(t)$ of a generic frontal curve $f : I \to \R^3$ 
appears just on $\{ s = 0 \}$. The list of singularities on 
$(t, 0; r) \in I\times \{ 0\} \times \R$ is given by, up to diffeomorphisms,  
$$
\begin{array}{|c|c|c|}
\hline
 {\mbox{\rm type}}  & {\mbox{\rm singularities}} & {\mbox{\rm codimension}}
\\
\hline
(1, 2, 3) & {\mbox{\rm cuspidal edge (CE$_{2,3}$)}} & 1
\\
(1, 2, 4) & {\mbox{\rm folded umbrella (FU$_{2,3}$)}} &  2
\\
(2, 3, 4) & {\mbox{\rm swallowtail (SW$_{2,3}$)}} & 2
\\
(2, 3, 5) & {\mbox{\rm folded pleat (FP$_{2,3}$)}} & 3
\\
(3, 4, 5) & {\mbox{\rm cuspidal swallowtail (CSW$_{2,3}$)}} & 3
\\
\hline
\end{array}
$$
\ent

The first column indicates the possible types of the directrix. 
The second column shows the names of right-left diffeomorphism class of singularities of parallels. 
The third column gives codimension of the locus of the type on the $(t, s; r)$-space $I\times\R\times\R$.

{\footnotesize 

  \begin{center}
%
\includegraphics[width=2.5truecm, height=2truecm, clip, 
]{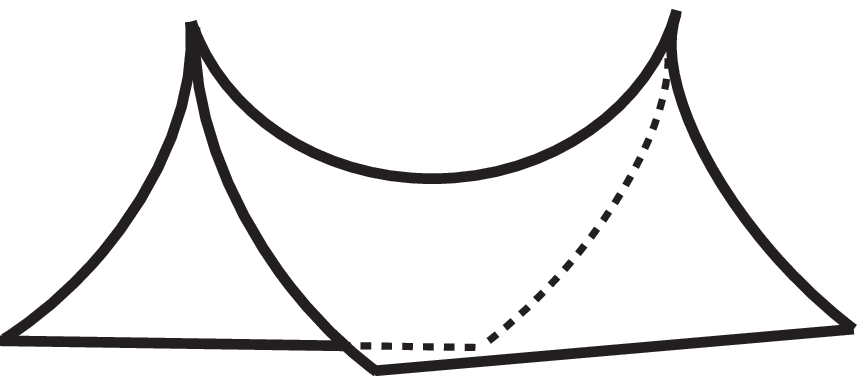} 
      \hspace{0.3truecm}
            \includegraphics[width=2.5truecm, height=2truecm, clip, 
]{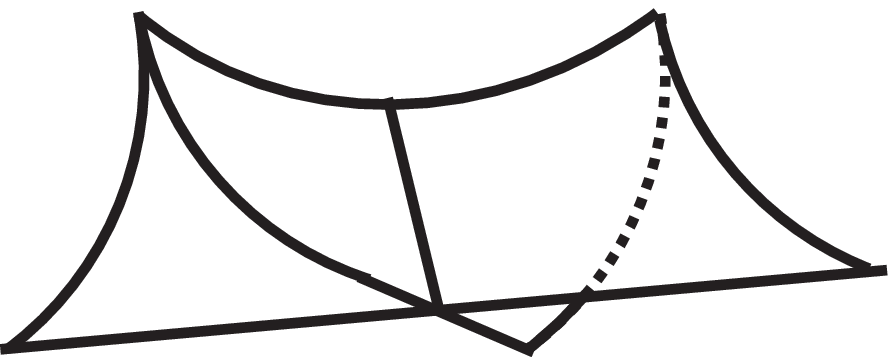} 
              \hspace{0.3truecm}
            \includegraphics[width=2.5truecm, height=2truecm, clip, 
            ]{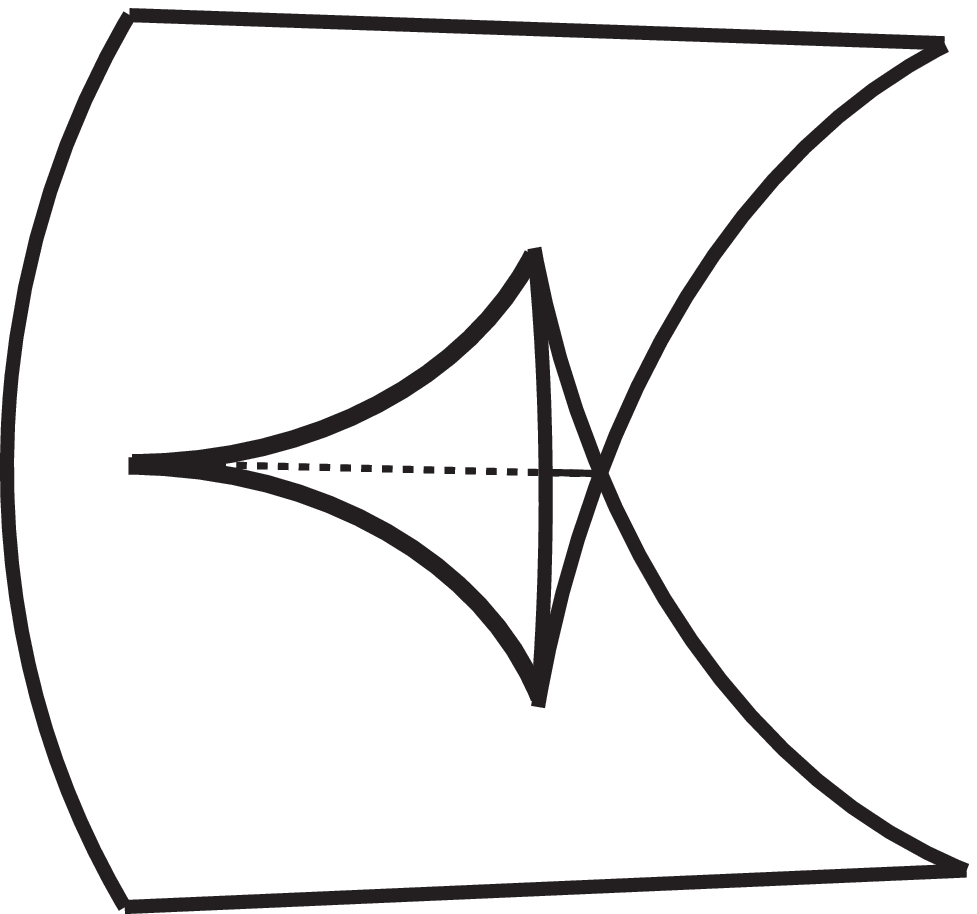} 
            \hspace{0.2truecm}
                  \includegraphics[width=2.5truecm, height=2truecm, clip, 
]{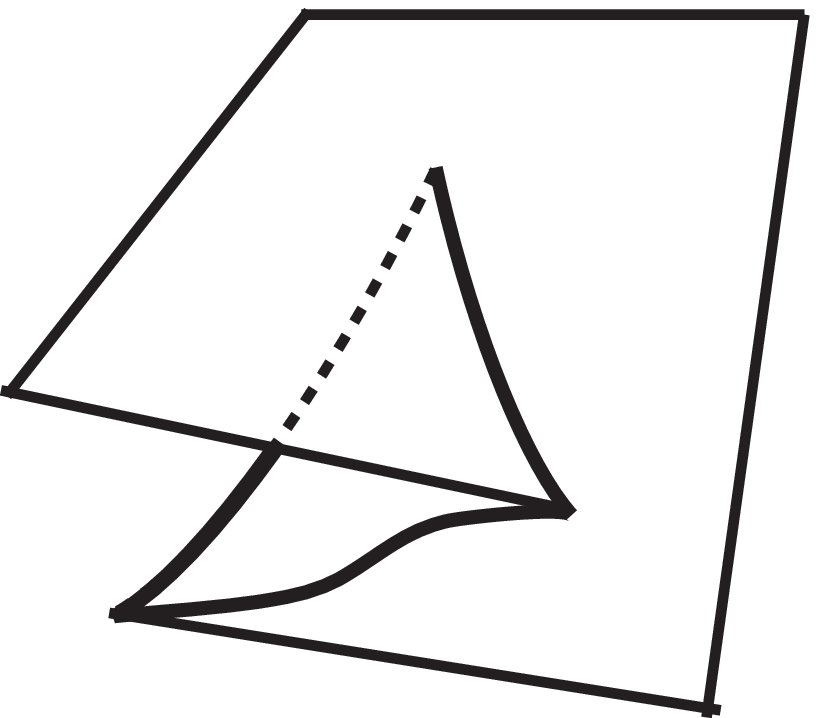}
         \hspace{0.3truecm}
 \includegraphics[width=2.5truecm, height=2truecm, clip, 
]{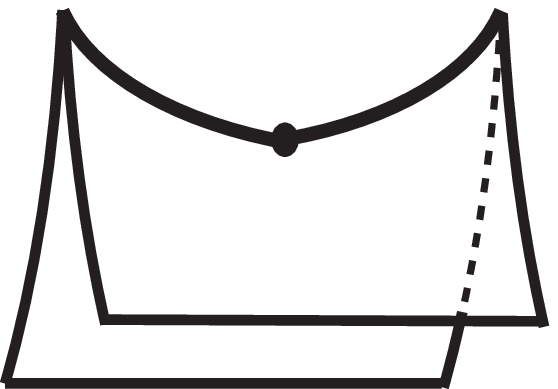} 
            \\
          \hspace{0.3truecm}
      cuspidal edge   
      \hspace{0.9truecm}
      folded umbrella 
        \hspace{0.9truecm}
        swallowtail
       \hspace{1truecm}
    folded pleat 
     \hspace{0.7truecm}
     cuspidal swallowtail
      \end{center}
    
}

It is known that the diffeomorphism type of tangent developables to curves in $\R^3$ 
of types in Theorem \ref{classification-p=2} is uniquely determined except for the case
$(2, 3, 5)$ (\cite{Ishikawa95}). 
The tangent surfaces of curves of type $(2, 3, 5)$ fall into two diffeomorphism classes, the generic 
folded pleat and the non-generic folded pleat (\cite{IMT}). 
We call a map-germ {\it the generic folded pleat} if it is diffeomorphic (right-left equivalent) to 
a map-germ $(\R^2, 0) \to (\R^3, 0)$ defined by 
$$
(u, t) \mapsto (u, \ \ t^3 + ut + \frac{3}{4}t^4 + \frac{1}{2}ut^2, \ \ \frac{3}{5}t^5 + \frac{1}{3}ut^3 + \frac{1}{2}t^6 + \frac{1}{4}ut^4). 
$$
Note that the above normal form is diffeomorphic to that given in \cite{IMT}. 

\ber
{\rm
The folded pleat and the cuspidal swallowtail appear at an isolated point on the directrix instantaneously. 
There is a bifurcation of a folded pleat (resp. a cuspidal swallowtail) depicted as follows: 
\begin{center}
\includegraphics[width=15truecm, height=3truecm, clip, 
]{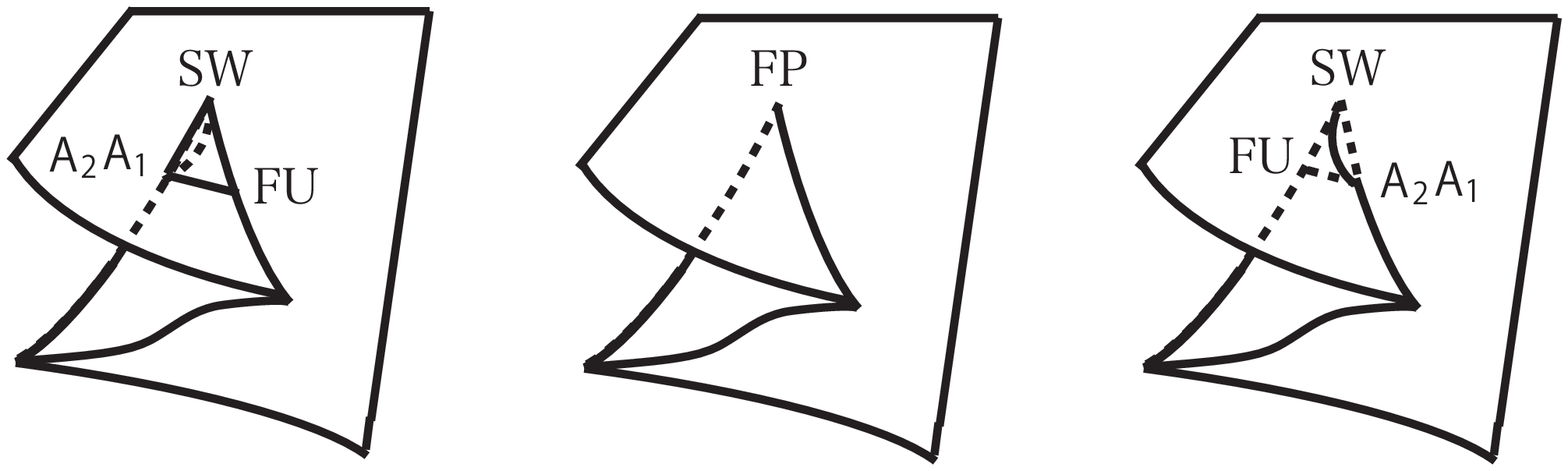} 
\\
A bifurcation of a folded pleat. 
\end{center}
\vspace{-0.5truecm}
\begin{center}
\includegraphics[width=16truecm, height=3truecm, clip, 
]{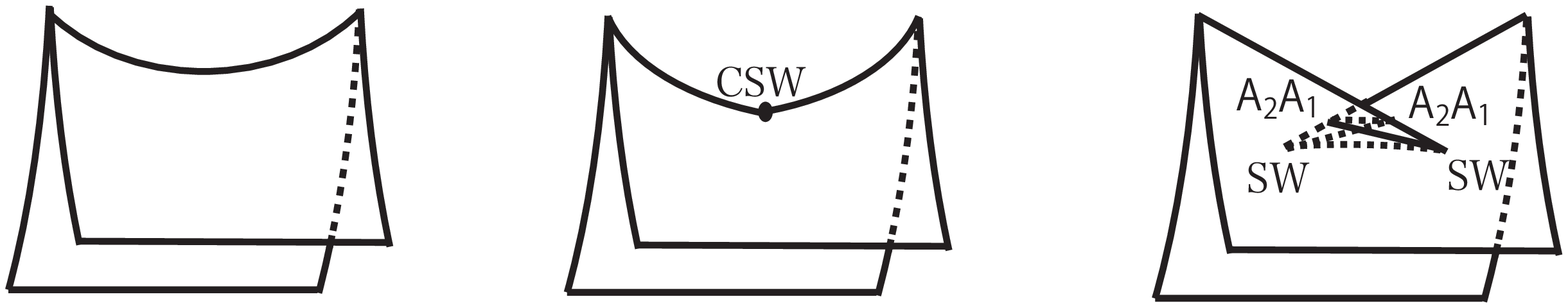} 
\\
A bifurcation of cuspidal swallowtail. 
\end{center}

Here the folded pleat bifurcates to a swallowtail, a folded umbrella and a transversal intersection of a cuspidal edge with a regular surface (denoted by $A_2A_1$). The cuspidal swallowtail bifurcates to two swallowtails and two $A_2A_1$ or to just a cuspidal edge. 
}
\enr

\

Let $p = 3$. Then the possible types are 
$$
(1, 2, 3, 4), (1, 2, 3, 5), (2, 3, 4, 5), (2, 3, 4, 6), (3, 4, 5, 6)
$$ 
by Theorem \ref{classification-curve-type}. 
Thus we have

\bet
\label{classification-p=3} 
{\rm (Classification of singularities of generic parallels to tangent maps of frontal curves in $\R^4$.)} \ 
Let $p = 3$. 
Any singular point on parallels
$f(t) + s\tau(t) + r\nu(t)$ 
along a generic normally parallel normal $\nu(t)$ of 
the tangent surface $f(t) + s\tau(t)$ of a generic frontal curve $f : I \to \R^4$ 
appears just on $\{ s = 0\}$. The list of singularities on 
$(t, 0 ; r) \in I\times \{ 0\}\times\R$ is given by 
$$
\begin{array}{|c|c|c|}
\hline
{\mbox{\rm type}} & {\mbox{\rm singularity}}  & {\mbox{\rm codimension}}
\\
\hline
(1, 2, 3, 4)  & {\mbox{\rm cuspidal edge (CE$_{2,4}$)}} & 1
\\
(1, 2, 3, 5) & {\mbox{\rm cuspidal edge (CE$_{2,4}$)}} &  2
\\
(2, 3, 4, 5) & {\mbox{\rm open swallowtail (OSW$_{2,4}$)}} & 2
\\
(2, 3, 4, 6) & {\mbox{\rm  unfurled swallowtail (USW$_{2,4}$)}} & 3
\\
(3, 4, 5, 6) & {\mbox{\rm cuspidal swallowtail (CSW$_{2,4}$)}} &  3
\\
\hline
\end{array}
$$
\ent

It is known the diffeomorphism type of tangent surfaces to curves 
of types $(1, 2, 3, 4), (1, 2, 3, 5)$ and $(2, 3, 4, 5)$ in $\R^4$ 
is uniquely determined (see \cite{Ishikawa12}). 
We show that it is also the case for $(3, 4, 5, 6)$ (Proposition \ref{cuspidal-swallowtail-normal-form} 
of \S \ref{The cuspidal swallowtail surfaces}). 

The singularities of  tangent surfaces to curves of type $(2, 3, 4, 6)$ 
turn to fall into several diffeomorphism classes (see \S \ref{The swallowtail and its openings}). 
The generic normal form is given by 
$$
(t, u) \mapsto 
(u, \ \ t^3+ut, \ \ \frac{3}{4}t^4+\frac{1}{2}ut^2, \ \ \frac{3}{7}t^7+\frac{1}{5}ut^5), 
$$
which is called the {\it unfurled swallowtail}.

\begin{center}

\includegraphics[width=2.2truecm, height=2truecm, clip, 
]{cuspidal-edge.eps} 
\hspace{0.7truecm}
\includegraphics[width=2.2truecm, height=2truecm, clip, 
]{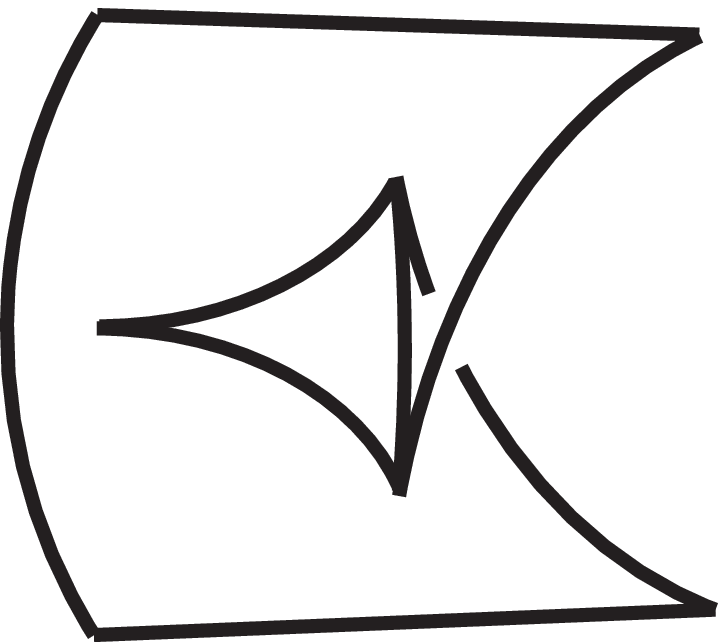} 
\hspace{0.7truecm}
\includegraphics[width=2.2truecm, height=2truecm, clip, 
]{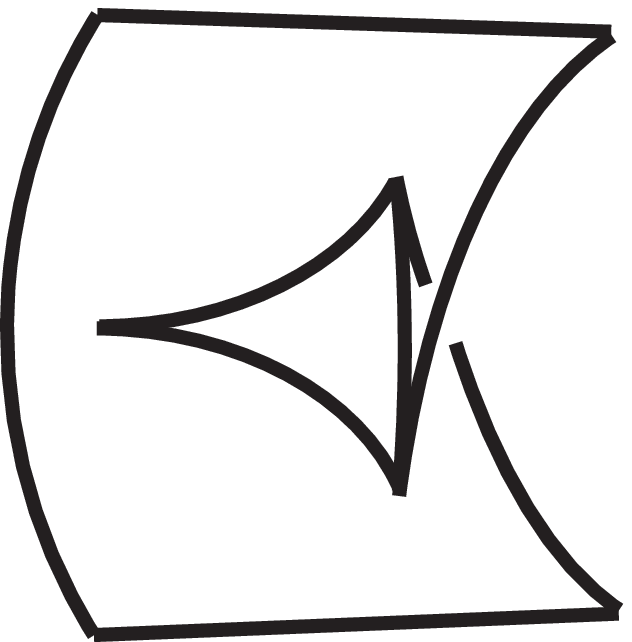}
\hspace{0.7truecm}
\includegraphics[width=2.2truecm, height=2truecm, clip, 
]{CS.eps}
\\
{\footnotesize 
CE \hspace{2truecm} OSW \hspace{2truecm}USW \hspace{2truecm} CSW
}

\end{center}

The following bifurcation occurs by parallel deformations of tangent surfaces of curves of type 
$(2, 3, 5, 6)$ (resp. $(3, 4, 5, 6)$): 

\begin{center}
\includegraphics[width=15truecm, height=3truecm, clip, 
]{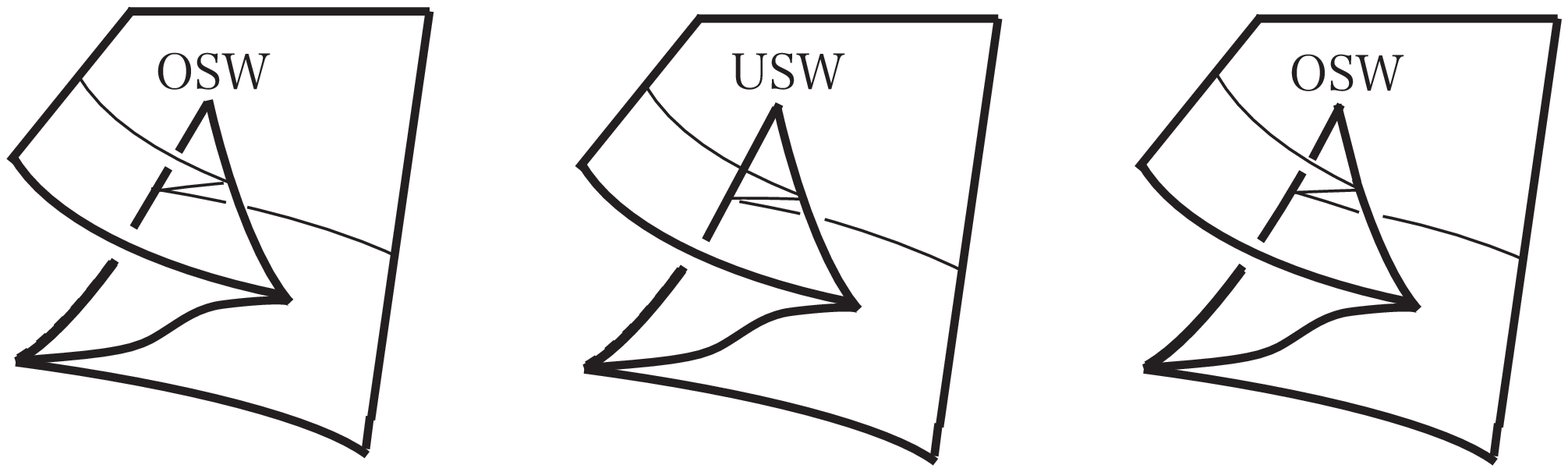}
\\
Bifurcation of the unfurled swallowtail. 
\end{center}
\vspace{-1truecm}
\begin{center}
\includegraphics[width=17truecm, height=3truecm, clip, 
]{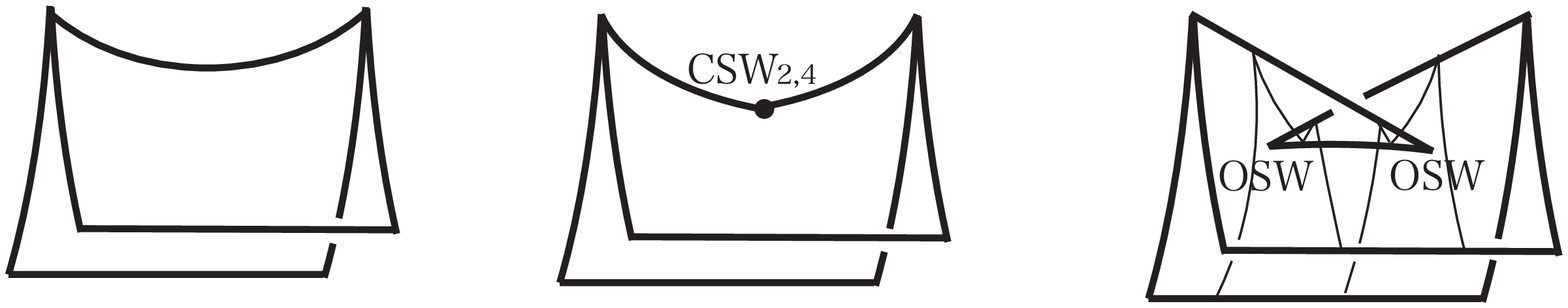} 
\\
Bifurcation of the cuspidal swallowtail in $\R^4$. 
\end{center}

\section{Parallels to tangent surfaces of frontal curves} 
\label{Parallels to tangent surfaces of frontal curves}

We recall here a sufficient condition on a curve to be a frontal. 

\bel
\label{sufficient-condition-frontal-tangent-frontal}
Let $f : (\R^1, 0) \to \R^{1+p}$ be a curve-germ. 
Suppose that 
$$
a_1 := \min \{ k \mid \rank(W_k(f)(0)) = 1\} < \infty, \ 
a_2 := \min \{ k \mid \rank(W_k(f)(0)) = 2\} < \infty. 
$$ 
Then we have the followings. 

{\rm (1)} 
$f$ is frontal and its tangent developable $\Tan(f) : (\R, 0)\times \R \to \R^{1+p}$ is a proper frontal.  
The singular locus of $\Tan(f)$ is equal to $(\R, 0)\times \{ 0\}$ if $a_2 - a_1 = 1$ and 
is equal to $(\R, 0)\times\{ 0\} \cup \{ 0\}\times\R$ if $a_2 - a_1 \geq 2$ 
{\rm (Lemma 2.15 of \cite{Ishikawa20b})}. 

{\rm (2)}
Let $\tau$ be a unit frame of $T_f$ and take a unit frame $\mu$ of $T_{\Tan(f)}\vert_{(\R, 0)\times 0}$ 
such that $\tau, \mu$ form an orthonormal frame of $T_{\Tan(f)}\vert_{(\R, 0)\times 0}$. 
Extend $\tau, \mu$ to an orthonormal frame of $T_{\Tan(f)}$ trivially by 
the constancy of $T_{\Tan(f)}$ along $t\times\R$ for any $t \in (\R, 0)$. 
Then $\tau'(t) = \kappa(t)\mu(t)$ for unique function $\kappa(t)$ on $(\R, 0)$ and 
$\ord_0\, \kappa = a_2 - a_1 - 1$. The curve $f$ has an inflection point at $t = 0$ 
if and only if $a_2 - a_1 \geq 2$. 
\enl

\Proof
(1) 
Note that $1 \leq a_1 < a_2$. 
There exist a $C^\infty$ coordinate $t$ of $(\R, 0)$ and an affine coordinates 
$x_1, x_2, \dots, x_{1+p}$ such that $f$ is given by 
$$
(x_1\circ f)(t) = t^{a_1}, \ \ (x_2\circ f)(t) = t^{a_2} + o(t^{a_2}), \ \ (x_i\circ f)(t) = o(t^{a_2}) \ (i = 3, \dots, 1+p). 
$$
Then we have a tangential field $\tau(t)$ for $t \not= 0$ by setting $\tau(t) = 
\frac{1}{a_1t^{a_1-1}}f'(t)$, which extends uniquely to 
a $C^\infty$ tangential field $\tau : (\R, 0) \to \R^{1+p}$ of form 
$$
\tau(t) = 
\left(	1, (a_2/a_1)t^{a_2-a_1} + o(t^{a_2-a_1}), \tau_3(t), 
\dots, \tau_{1+p}(t)\right), 
$$
with $\tau_i(t) = o(t^{a_2-a_1}), (i = 3, \dots, 1+p)$. Therefore $f$ is a frontal. 
Then the tangent surface of $f$ is given by $\Tan(f)(t, s) = f(t) + s\tau(t)$. 
We have that the Jacobian ideal of $\Tan(f)$ is principal and, in fact, 
is generated by $st^{a_2-a_1-1}$. Therefore $\Tan(f)$ has a nowhere dense singular locus 
and it turns out to be a proper frontal (see \cite{Ishikawa20}, Lemma 2.3). 

(2) We write $\tau'(t) = a(t)\tau(t) + \kappa(t)\mu(t)$ for some functions $a, \kappa$. 
Since $\tau'(t)\cdot\tau(t) = 0$ and $\mu(t)\cdot\tau(t) = 0$, we see $a(t) = 0$. 
The tangent map $\Tan(f)$ is given by $(t, s) \mapsto f(t) + s\tau(t)$ up to 
diffeomorphisms preserving $(\R, 0)\times 0$. The Jacobi matrix of $\Tan(f)$ 
is given by $(f'(t) + s\tau'(t), \tau(t))$. The ideal generated by $2\times 2$-minors of $(f'(t) + s\tau'(t), \tau(t))$ 
is equal to that for $s\kappa(t)(\mu(t), \tau(t))$ and it is generated by the function $s\kappa(t)$. 
By the argument in (1), we see that the ideals generated by $st^{a_2-a_1-1}$ and by $s\kappa(t)$. 
Therefore in particular we have that $\ord_0\, \kappa = a_2 - a_1 - 1$. Moreover 
$\kappa(0) = 0$ if and only if $a_2 - a_1 \geq 2$ as required. 
\qed

\

The types appeared in Theorem \ref{classification-curve-type} satisfy $a_2 - a_1 = 1$. 
Therefore we see that generically a frontal curve $f$ has no inflection points 
and that $\Tan(f)$ is a proper frontal 
by Lemma \ref{sufficient-condition-frontal-tangent-frontal}.

\bef
{\rm A curve $\tau : I \to \R^{1+p}$ is said to be of {\it primitive type} 
$(b_1, b_2, \dots, b_p, b_{1+p})$ if 
$b_i = \min\{ k \mid \rank(\widetilde{W}_k(\tau)(t_0)) = i\}, 1 \leq i \leq 1+p$, where $\widetilde{W}_k(\tau)(t) := 
(\tau(t), \tau'(t), \dots, \tau^{(k-1)}(t))$, $(1+p)\times k$-matrix, i.e. Wronskian of the primitive of $\tau$. 
}
\enf

\bel
\label{Primitive-type-and-type}
Let $f : I \to \R^{1+p}$ be a frontal curve and $\tau : I \to \R^{1+p}$ a tangential frame of $f$ 
with $f'(t) = a(t)\tau(t)$ for unique function $m(t)$. 
If $\tau$ is of primitive type $(b_1, b_2, \dots, b_p, b_{1+p})$ at $t = t_0$, 
and $\ord_{t_0}a = m$ i.e. $\tau(t_0) = \cdots = \tau^{(m-1)}(t_0) = 0, \tau^{(m)}(t_0) \not= 0$, then 
$f$ is of type $(m+b_1, m+b_2, \dots, m+b_p, m+b_{1+p})$ at $t = t_0$. 
\enl

\Proof
From $f'(t) = a(t)\tau(t)$, we have 
$f''(t) = a'(t)\tau(t) + a(t)\tau'(t)$, $f'''(t) = a''(t)\tau(t) + 2a'(t)\tau'(t) + a(t)\tau''(t)$ and so on. 
If $\ord_{t_0}a = m$, then $W_r(f)(t_0) = O$ and  
we have $\left( f^{(m+1)}(t_0), \dots, f^{(m+k)}(t_0)\right) =  \left( \tau(t_0), \dots, \tau^{(k)}(t_0)\right)\cdot A(t_0)$,  where 
$$
A = 
\left(
\begin{array}{ccccc}
a^{(m)} & a^{(m+1)} &  \cdots & \cdots & a^{(m+k-1)} 
\\
0 & a^{(m)}  & a^{(m+1)} & \cdots & \cdots
\vspace{0.2truecm}
\\
\vdots & \vdots & \ddots & \ddots & \vdots
\vspace{0.2truecm}
\\
0 & 0 & \cdots & a^{(m)} & a^{(m+1)}
\\
0 & 0 & \cdots & 0 & a^{(m)}
\end{array}
\right), 
$$
a $k\times k$-matrix. 
Then $W_m(f)(t_0) = O$ and $\rank\  W_{m+k}(f)(t_0) = \rank\  \widetilde{W}_k(\tau)(t_0)$. 
Therefore we have that $f$ is of type $(a_1, a_2, \dots, a_{p+1}) = (b_1 + m, b_2 + m, \dots, b_{p+1} + m)$ at $t_0$.  
\QED

\

\noindent
{\it Proof of Proposition \ref{generic-frontal-curves}:} 
Write $f'(t) = a(t)\tau(t)$ for a function $a : I \to \R$ and $\tau : I \to S^p \subset \R^{1+p}$. 
Note that $f$ is obtained by integration from $a$ and $\tau$, and $f$ is determined up to constant. 
The spherical curve $\tau$ is regarded as a curve in the projective space $\R P^p$ as well and 
the type of $\tau$ at a point $t \in I$ as an affine curve is equal to $(c_1, c_2, \dots, c_p)$ if and only if 
the primitive type of $\tau$ at $t$ as 
a curve in $\R^{1+p}$ is equal to $(1, 1+c_1, 1+c_2, \dots, 1+c_p)$. 
On the other hand it is known that generically a curve in $\R P^p$ is of type $(1, 2, \dots, p-1, p)$ 
or $(1, 2, \dots, p-1, p+1)$ and the latter occurs only at an isolated point
 (see \cite{Scherbak} and Corollary 5.2 of \cite{Ishikawa12}). 
 Note that the genericity condition is given by a transversality of jet sections to a closed semi-algebraic set 
 in the jet space $J^{p+1}(I, \R P^p)$ of codimension $\geq 2$. 
Thus by taking $\tau$ generically we have that there exist discrete set $D \subset I$ such that 
$\tau$ is of primitive type $(1. 2. \dots, p, p+1)$ at $t \in I \setminus D$ and 
$(1, 2, \dots, p, p+2)$ at $t \in D$. 
Then by taking the function $a$ generically we have that $a$ does not vanish at $D$ and 
any zero $t$ of $a$ on $I \setminus D$ is simple, i.e. $\ord_t a = 1$. 
Therefore, by Lemma \ref{Primitive-type-and-type}, we have generically that $f$ is of type 
$(1. 2. \dots, p, 1+p)$, $(1, 2, \dots, p, 2+p)$ or $(2, 3, \dots, 1+p, 2+p)$. 
\QED

\

\noindent
{\it Proof of Theorem \ref{parallel-of-tangent}:}
Let $f : (\R, 0) \to \R^{1+p}$ be a frontal curve and $F = \Tan(f) : 
(\R, 0) \times \R \to \R^{1+p}$ the tangent developable of $f$. 
Let $\{ \tau, \mu, \nu_1, \dots, \nu_{p-1}\}$ be an orthonormal frame of $F^*T\R^{1+p}$ along $f$. 
Here we demand that $\tau$ is a unit section of $T_f$, $\{ \tau, \mu\}$ is an orthonormal frame 
of $T_F\vert_{(\R, 0)\times 0}$ and $\{ \nu_1, \dots, \nu_{p-1}\}$ is a parallel orthonormal 
frame of $N_F\vert_{(\R, 0)\times 0}$.  Then we extend 
$\{ \tau, \mu, \nu_1, \dots, \nu_{p-1}\}$ trivially to the frame of $F^*T\R^{1+p}$ 
by the constancy of $T_F\vert_{t \times \R}$. 
Then we have $f'(t) = a(t) \tau(t)$ and the structure equation of the moving frame: 
$$
\left(
\begin{array}{c}
\tau'(t) \\
\mu'(t) \\
\nu'_1(t) \\
\vdots \\
\nu'_{p-1}(t)
\end{array}
\right)
= 
\left(
\begin{array}{ccccc}
0 & \kappa(t) & 0 & \cdots & 0 \\
- \kappa(t) & 0 & \ell_1(t) & \cdots & \ell_{p-1}(t) \\
0 & -\ell_1(t) & 0 & \cdots & 0 \\
\vdots & \vdots & \vdots & \ddots & \vdots \\
0 & -\ell_{p-1}(t) & 0 & \cdots & 0
\end{array}
\right)
\left(
\begin{array}{c}
\tau(t) \\
\mu(t) \\
\nu_1(t) \\
\vdots \\
\nu_{p-1}(t)
\end{array}
\right)
$$ 
for uniquely determined functions $a(t), k(t), \ell_1(t), \dots, \ell_{p-1}(t)$. 

Let $\nu$ be a normally parallel, not necessarily unit, normal field along $\Tan(f)$. Then 
$\nu = \sum_{i=1}^{p-1} r_i\nu_i$ for some $r_i \in \R\ (i = 1, \dots, p-1)$. 
Consider the parallel 
$$
P(t, s) := {\mbox{\rm{P}}}_\nu(F)(t, s) = F + \nu = f(t) + s\tau(t) + \nu(t), 
$$ 
of $F = \Tan(f)$ generated by $\nu$. 
Then the Jacobi matrix of $P$ is given by 
$$
\left(f'(t) + s\tau'(t) + \sum_{i=1}^{p-1} r_i\nu_i'(t), \ \tau(t)\right) 
= \left(a(t)\tau(t) + \left\{s\kappa(t) - \sum_{i=1}^{p-1} r_i\ell_i(t)\right\}\mu(t), \ \tau(t)\right), 
$$
the rank of which is equal to that of the matrix
$\left(\left\{s\kappa(t) - \sum_{i=1}^{p-1} r_i\ell_i(t)\right\}\mu(t), \ \tau(t)\right)$. 
Therefore the singular locus $\Sigma(P)$ of the parallel $P$ is given by 
$\{ (t, s) \in (\R, 0)\times \R \mid s \kappa(t) - \sum_{i=1}^{p-1}r_i\ell_i(t) = 0\}$. 
By the assumption, $\kappa(t)$ does not vanish. Therefore  
$\Sigma(P)$ is parametrised by $s = \sum_{i=1}^{p-1}r_i\ell_i(t)/\kappa(t)$. 
Set 
$$
g(t) := f(t) + \sum_{i=1}^{p-1}r_i\frac{\ell_i(t)}{\kappa(t)}\tau(t) + \nu(t).
$$
Then 
$$ 
\begin{array}{rcl}
g'(t) & = & f'(t) + \sum_{i=1}^{p-1}r_i\left(\frac{\ell_i(t)}{\kappa(t)}\tau(t)\right)' + \nu'(t) 
\\
& = & a(t)\tau(t) + \sum_{i=1}^{p-1}r_i\left(\frac{\ell_i(t)}{\kappa(t)}\right)'\tau(t) 
+ \sum_{i=1}^{p-1}r_i\left(\frac{\ell_i(t)}{\kappa(t)}\right)\kappa(t)\mu(t) + \sum_{i=1}^{p-1}r_i(-\ell_i(t)\mu(t))
\\
& = & \left( a(t) + \sum_{i=1}^{p-1} r_i\left\{ \frac{\ell_i(t)}{\kappa(t)}\right\}'\right)\tau(t).
\end{array}
$$
Thus $g$ has the same tangent frame $\tau$ with $f$. The tangent developable of $g$ is given by 
$$
\Tan(g)(t, s) = f(t) + \sum_{i=1}^{p-2}r_i\frac{\ell_i(t)}{\kappa(t)}\tau(t) + \nu(t) + s\tau(t) 
= f(t) + \left(s + \sum_{i=1}^{p-2}r_i\frac{\ell_i(t)}{\kappa(t)}\right)\tau(t) + \nu(t). 
$$
Then by the diffeomorphism $(t. s) \mapsto (t, s + \sum_{i=1}^{p-2}r_i\frac{\ell_i(t)}{\kappa(t)})$ 
on $(\R, 0)\times \R$, we have that ${\mbox{\rm{P}}}_\nu(\Tan(f))$ is right equivalent 
to $\Tan(g)$. 
\qed

\

\noindent
{\it Proof of Theorem \ref{classification-curve-type}:} 
By Theorem \ref{parallel-of-tangent}, we have
$g'(t) = \left( a(t) + \sum_{i=1}^{p-1} r_i\left\{ \frac{\ell_i(t)}{\kappa(t)}\right\}'\right)\tau(t)$. 
We set $b(t) = a(t) + \sum_{i=1}^{p-1} r_i\left\{ \frac{\ell_i(t)}{\kappa(t)}\right\}'$. 
Let us consider the system of equations:
$$
\left\{
\begin{array}{ccccccccc}
a_0(t) & + & a_1(t)\, r_1 & +  & \cdots & + & a_{p-1}(t)\, r_{p-1} & = & 0, 
\\
a_0'(t) & + & a_1'(t)\, r_1 & +  & \cdots  & + & a_{p-1}'(t)\, r_{p-1} & = & 0, 
\\
a_0''(t) & + & a_1''(t)\, r_1 & +  & \cdots & + & a_{p-1}''(t)\, r_{p-1} & = & 0, 
\end{array}
\right.
$$
for a curve $\alpha : I \to \R^p, \alpha(t) = (a_0(t), a_1(t), \dots, a_{p-1}(t))$ and 
$r = (r_1, r_2, \dots, r_{p-1}) \in \R^{p-1}$. 
By the transversality theorem, the curves $\alpha$ satisfying the following properties form an open dense set: 
The set of $r$ satisfying the above system of equations for some $t \in I$ 
form a semi-algebraic subset $\Sigma$ in $I\times \R^{p-1}$ 
of dimension $\leq p - 3$ and $a_0(t) \not= 0$ for any $t \in I$. 
Then $\pi(\Sigma) \subset \R^{p-1}$ is a semi-algebraic subset of dimension $\leq p-3$ in $\R^{p-1} \setminus \{ 0\}$ which does not pass through the origin, where $\pi : I \times \R^{p-1}, \pi(t, r) = r$. 
Then $\Pi(\pi(\Sigma)) \subset \R P^{p-2}$ is a semi-algebraic subset of dimension $\leq p-3$ in $\R P^{p-2}$. 

Now let $f : I \to \R^{1+p}$ be any frontal. Set $f' = a\tau$ for a function $a : I \to \R$ and 
$\tau : I \to S^p \subset \R^{1+p}$. Then perturb $\tau$ to a generic map as in Proposition \ref{generic-frontal-curves}. 
Let $\tau(t), \mu(t), \nu_1(t), \dots, \nu_{p-1}(t)$ be the corresponding moving frame of $f$ and 
$\kappa(t), \ell_1(t), \dots, \ell_{p-1}(t)$ the system of invariants (curvature and torsions) of $f$. 
Perturb $a(t), \kappa(t), \ell_1(t), \dots, \ell_{p-1}(t)$ such that 
$a_0(t) = a(t), a_1(t) = (\ell_1(t)/\kappa(t))', \dots, a_{p-1}(t) = (\ell_{p-1}(t)/\kappa(t))'$ satisfies the genericity condition. 
Then we perturb the frame $\tau(t), \mu(t), \nu_1(t), \dots, \nu_{p-1}(t)$ by the structure equation. 
Finally $f(t)$ is perturbed by integration of $a(t)\tau(t)$ for the perturbed $a(t)$ and $\tau(t)$. 

Take any $R \in \R^{1+p}, R \not= 0$ with $[R] \in \R P^{p-2} \setminus \Pi(\pi(\Sigma))$ and set 
$\nu := \sum_{i=1}^{p-1} R_i\nu_i$. Then, for any parallel $P_{r\nu} = f + r\nu$, for any $r \in \R$, we have 
$P_{r\nu}'(t) = b(t)\tau(t)$ with 
$b(t) = a(t) + \sum_{i=1}^{p-1} rR_i(\ell_i(t)/\kappa(t))'$ 
such that $b(t)$ has order $\leq 2$ for any $t \in I$. 

We have that there exist discrete set $D \subset I$ such that 
$\tau$ is of primitive type $(1. 2. \dots, p, p+1)$ at $t \in I \setminus D$ and 
$(1, 2, \dots, p, p+2)$ at $t \in D$. Moreover, we can arrange $a(t)$ slightly so that 
the discrete point set $\{ (a(t), a'(t)) \mid t \in D\}$ does not intersect with the union of lines 
$\bigcup_{t \in D} \R(\sum_{i=1}^{p-1} R_i(\ell_i(t)/\kappa(t))', \sum_{i=1}^{p-1} R_i(\ell_i(t)/\kappa(t))'')$   
in $\R^2$. 
Then we have 
$(b(t), b'(t)) \not= (0, 0)$ for any $r \in \R$ and $t \in D$. 
Thus $b(t)$ has order $\leq 1$ for any $t \in I \setminus D$. 
Therefore, by Lemma \ref{Primitive-type-and-type}, we have the required result on the possible types appearing 
in directrixes. 
\QED

\section{The cuspidal swallowtail surfaces in $\R^4$} 
\label{The cuspidal swallowtail surfaces}

The determinacy of generating families of Legendre varieties associated to 
curves of type $(3, 4, 5, 6)$ is shown in \cite{Ishikawa93}. See also \cite{Ishikawa95}. 
Using it, we give the normal form of tangent surfaces to curves of type $(3, 4, 5, 6)$ 
in $\R^4$. 

\bep
\label{cuspidal-swallowtail-normal-form}
Tangent surfaces of curve-germs of type $(3, 4, 5, 6)$ in $\R^4$ has unique diffeomorphism type. 
In fact they are diffeomorphic to the germ $(\R^2, 0) \to \R^4$ defined by 
$$
{\mbox{\rm CSW}}_{2,4}: \ 
(u, t) \mapsto \left( u, \ \ t^4 + ut, \ \ \frac{4}{5}t^5 + \frac{1}{2}ut^2, \ \ \frac{2}{3}t^6 + \frac{1}{3}ut^3\right). 
$$
\enp

\

\noindent
{\it Proof of Proposition \ref{cuspidal-swallowtail-normal-form}.} 
By Theorem 1 of \cite{Ishikawa93}, it is known that the ruled $3$-fold by osculating planes 
to any curve of type $(3,4,5,6)$ in $\R^4$ is diffeomorphic to the envelope of 
the family of hyperplanes 
$$
t^6 + x_1 t^3 + x_2 t^2 + x_3 t + x_4 = 0, 
$$
with parameter $t$. Then the tangent surface of the curve is diffeomorphic to the singular locus of the envelope. 
The envelope is given by 
$(x_1, x_2, t) \mapsto (x_1, x_2, x_3, x_4)$, where 
$$
x_3 = - 6t^5 - 3x_1 t^2 - 2x_2 t, \quad x_4 = 5t^6 + 2x_1 t^3 + x_2 t^2. 
$$
The singular locus is given by $x_2 = - 15t^4 - 3x_1 t$ and therefore by 
the map-germ
$$
\begin{array}{rcl}
(x_1, t) & \mapsto & \left( x_1, \ - 15t^4 - 3x_1 t, 24t^5 + 3x_1 t^2, \ -10t^6 - x_1t^3 \right) 
\\
& = & \left( x_1, \ - 15(t^4 + \frac{1}{5}x_1 t), 30(\frac{4}{5}t^5 + \frac{1}{10}x_1 t^2), \ -15(\frac{2}{3}t^6 + \frac{1}{15}x_1t^3)\right), 
\end{array}
$$
which is diffeomorphic to the normal form by setting $u = \frac{1}{5}x_1$. 
\QED

\section{The swallowtail and its openings in $\R^4$} 
\label{The swallowtail and its openings} 

The normal form of the swallowtail surface is given by 
$$
{\mbox{\rm SW}}_{2,3} : \ (t, u) \mapsto \left( u, \ \ t^3 + ut, \ \ \frac{3}{4}t^4 + \frac{1}{2}ut^2\right). 
$$
The tangent surface to any curve of type $(2, 3, 4, 5)$ is diffeomorphic to the 
{\it open swallowtail}, which is the \lq\lq versal" opening of the swallowtail surface in $\R^3$. 
$$
{\mbox{\rm OSW}}_{2,4}: \
(t, u) \mapsto \left( u, \ \ t^3 + ut, \ \ \frac{3}{4}t^4 + \frac{1}{2}ut^2, \ \ \frac{3}{5}t^5 + \frac{1}{3}ut^3\right), 
$$
For the details and proofs, see \cite{Ishikawa12}. 

\

Related to Theorem \ref{classification-p=3}, we are lead to study the classification for more degenerate \lq\lq opening" 
of the swallowtail.

\bet
\label{unfurled-swallowtail-classification}
The tangent surface $\Tan(f)$ of a curve $f$ of type $(2, 3, 4, 6)$ 
is diffeomorphic to 
$$
(t, u) \mapsto 
\left( u, \ \ t^3+ut, \ \ \frac{3}{4}t^4+\frac{1}{2}ut^2, \ \ \psi(u, t^3+ ut)(\frac{3}{7}t^7+\frac{1}{5}ut^5)\right), 
$$
for some function $\psi : (\R^2, 0) \to \R$. 
If $\psi$ is identically zero, then the tangent surface $\Tan(f)$ is diffeomorphic to 
the embedded swallowtail 
$$
{\mbox{\rm SW}}_{2,4} : \ (t, u) \mapsto (u, \ \ t^3 + ut, \ \ \frac{3}{4}t^4 + \frac{1}{2}ut^2, \ \  0), 
$$
in $\R^4$.  If $\psi(0, 0) \not= 0$, then $\Tan(f)$ is diffeomorphic to 
$$
{\mbox{\rm USW}}_{2,4} : \ 
(t, u) \mapsto 
(u, \ \ t^3+ut, \ \ \frac{3}{4}t^4+\frac{1}{2}ut^2, \ \ \frac{3}{7}t^7+\frac{1}{5}ut^5), 
$$
\ent

\ber
{\rm
The three singularities ${\mbox{\rm SW}}_{2,4}, {\mbox{\rm OSW}}_{2,4}$ and ${\mbox{\rm USW}}_{2,4} : (\R^2, 0) 
\to \R^4$ are 
not diffeomorphic to each other. 
}
\enr

\

To show Theorem \ref{unfurled-swallowtail-classification}, let 
$g(t, u) = (u, t^3 + ut)$. 
Set 
$$
{\mathcal J}_g := \langle du, d(t^3+ ut)\rangle_{{\mathcal E}_2} 
= \langle du, (3t^2+u)dt\rangle_{{\mathcal E}_2}, 
\quad 
{\mathcal R}_g := \{ h \in {\mathcal E}_2 \mid dh \in {\mathcal J}_g\}
$$
(see \cite{Ishikawa18, Ishikawa20}) and set 
$$
{\mathcal R}_g^{(i+1)} := 
\{ h \in t^{i+1} {\mathcal E}_2 \mid dh \in t^i{\mathcal J}_g \}, \ (i = 0, 1, 2, \dots). 
$$
Here ${\mathcal E}_2$ means the $\R$-algebra which consists of all $C^\infty$ function-germs $(\R^2, 0) \to \R$. 
Then we have 
$$
{\mathcal E}_2 \supset {\mathcal R}_g^{(1)} \supset {\mathcal R}_g^{(2)} \supset \cdots . 
$$

Let $T_0 := t^3 + ut = \dint{}{}(3t^2+u)dt$ and 
$$
T_i := \dint{}{} t^i(3t^2+u)dt = \dint{}{}(3t^{i+2} + ut^i)dt = 
\frac{3}{i+3}t^{i+3} + \frac{1}{i+1}ut^{i+1}, 
$$ 
for $i = 0, 1, 2, 3, \dots$. 
Then $T_i \in {\mathcal R}_g^{(i+1)}$. 

\bel
We have 
$$
{\mathcal R}_g^{(i+1)} = \left\{ h \in t^{i+1}{\mathcal E}_2 \ \left\vert\  \frac{\pa h}{\pa t} \in t^i \frac{\pa T}{\pa t}{\mathcal E}_2 \right\}\right.
= 
{\mathcal R}_g \cap t^{i+1}{\mathcal E}_2. 
$$
\enl

\Proof
Two inclusions $\subseteq$ are clear. To show the converse inclusions, let 
$h \in {\mathcal R}_g \cap t^{i+1}{\mathcal E}_2$. 
Then $dh = a du + b(3t^2+u)dt$ for some $a, b \in {\mathcal E}_2$. 
Moreover $h = t^{i+1}k$ for some $k \in {\mathcal E}_2$. 
Then 
$dh = t^{i+1}dk + (i+1)t^i k dt = t^{i+1}\frac{\pa k}{\pa u}du + t^{i}(t\frac{\pa k}{\pa t} + (i+1)k)dt$. Then $a = t^{i+1}\frac{\pa k}{\pa u}$, 
$b(3t^2 + u) = t^{i}(t\frac{\pa k}{\pa t} + (i+1)k)$ and we have $b \in t^i{\mathcal E}_2$. 
Therefore we have $\frac{\pa h}{\pa t} \in t^i \frac{\pa T}{\pa t}{\mathcal E}_2$ and $dh \in t^i{\mathcal J}_g$. 
Thus we have two inclusions $\supseteq$. 
\QED

\

By Lemma 2.4 (1') of \cite{Ishikawa95}, we have 

\bel
\label{algebraic-lemma}
{\rm
The ${\mathcal E}_2$-module ${\mathcal R}_g^{(4)}$ is generated by 
$T_3, T_4, T_5$ over $g^* : {\mathcal E}_2 \to {\mathcal E}_2$. 
}
\enl

By direct calculations, we have 

\bel
\label{algebraic-lemma-2}
{\rm 
$$
T_{2j+1} = \frac{j+2}{6} T_{j-1}^2 - \frac{j+2}{3j} u T_{2j-1}, \quad (j = 1, 2, 3, \dots). 
$$
In particular we have 
$$
T_3 = \frac{1}{2} T_0^2 - uT_1, \quad T_5 = \frac{2}{3} T_1^2 - \frac{2}{3} uT_3 = - \frac{1}{3}uT_0^2 + \frac{2}{3}u^2T_0 + \frac{2}{3}T_1^2. 
$$
}
\enl

\noindent
{\it Proof of Theorem \ref{unfurled-swallowtail-classification}:} 
Let $f : I \to \R^4$ be of type $(2, 3, 4, 6)$ at $t_0 \in I$. 
Then there exist a coordinate $t$ of $I$ centred at $t_0$ and a system of affine coordinates on $\R^4$ centred at $f(t_0)$ 
such that $f$ is expressed as
$$
f : x_1 = t^2, \ \ x_2 = t^3 + \varphi_2(t), \ \ x_3 = t^4 + \varphi_3(t), \ \ x_4 = t^6 + \varphi_4(t), 
$$
with $\ord_0 \varphi_2 > 3, \ord_0 \varphi_3 > 4, \ord_0 \varphi_4 > 6$. 
Then a tangential frame of $f$ is given by 
$$
\tau(t) =  
\left( 1, \ \ \frac{3}{2}t + \psi_2(t), \ \ 2t^2 + \psi_3(t), \ \ 3t^4 + \psi_4(t)\right), 
$$
with $\varphi_2'(t) = 2t\psi_2(t), \varphi_3'(t) = 2t\psi_3(t), \varphi_4'(t) = 2t\psi_4(t)$ and 
$\ord_0 \psi_2 > 1, \ord_0 \psi_3 > 2, \ord_0 \psi_4 > 4$. 
Then $\Tan(f) = f + s\tau$ is given by 
$$
x_1 = t^2 + s, x_2 = t^3 + \frac{3}{2}st + \varphi_2 + s\psi_2, x_3 = t^4 + 2st^2 + \varphi_3 + s\psi_3, 
x_4 = t^6 + 3st^4 + \varphi_4 + s\psi_4. 
$$
Set $u = -3(t^3+3)$. Then we have $x_1 = -3u$ and 
$$
\begin{array}{rcl}
x_2 & = & - \frac{1}{2}t^3 - \frac{1}{2}ut + \varphi_2 - \frac{1}{3}(3t^2 + u)\psi_2 = - \frac{1}{2}(T + \Phi_2(t, u)), 
\\
x_3 & = & - t^4 - \frac{2}{3}ut^2 + \varphi_3 - \frac{1}{3}(3t^2 + u)\psi_3 = -\frac{4}{3}(T_1 + \Phi_3(t, u)), 
\\
x_4 & = & - 2t^6 - ut^4 + \varphi_4 - \frac{1}{3}(3t^2 + u)\psi_4 = -4(T_3 + \Phi_4(t, u)), 
\end{array}
$$
with $\Phi_2 \in {\mathcal R}_g^{(2)}, \Phi_3 \in {\mathcal R}_g^{(3)}$ and $\Phi_4 \in {\mathcal R}_g^{(5)}$. 
Note that $T_0 \in {\mathcal R}_g^{(1)}, T_1 \in {\mathcal R}_g^{(2)}$ and $T_3 \in {\mathcal R}_g^{(4)}$. 
The germ of $\Tan(f)$ at $0$ is diffeomorphic to $\widetilde{F} : (\R^2, 0) \to (\R^4, 0)$ defined by 
$$
x_1 = u, \ \ x_2 = T_0 + \Phi_2(t, u), \ \ x_3 = T_1 + \Phi_3(t, u), \ \ x_4 = T_3 + \Phi_4(t, u), 
$$
by the affine transformation $(x_1, x_2, x_3, x_4) \mapsto (\frac{1}{3}x_1, -2x_2, -\frac{3}{4}x_3, -\frac{1}{4}x_4)$ on $\R^4$. 

Denote by $\widetilde{g}$ (resp. $\widetilde{G}$) the map-germ $(\R^2, 0) \to (\R^2, 0)$ (resp. 
$(\R^2, 0) \to (\R^3, 0)$) defined by $(t, u) \mapsto (u, T_0 + \Phi_2)$ (resp. $(t, u) \mapsto (u, T_0 + \Phi_2, T_1 + \Phi_3(t, u))$). 
Then, as a geometric meaning, $\widetilde{g}$ (resp. $\widetilde{G}$, $\widetilde{F}$) 
corresponds to tangent map of curves of type $(2, 3)$ (resp. $(2, 3, 4)$, $(2, 3, 4, 6)$). 
Moreover, as an algebraic aspect, we see that ${\mathcal J}_{\widetilde{g}} = {\mathcal J}_{\widetilde{G}} = 
{\mathcal J}_{\widetilde{F}} = {\mathcal J}_g$ and therefore 
${\mathcal R}_{\widetilde{g}} = {\mathcal R}_{\widetilde{G}} = 
{\mathcal R}_{\widetilde{F}} = {\mathcal R}_g$. It is known that the tangent surfaces of curves of type $(2, 3, 4)$ 
are diffeomorphic to the swallowtail surface. In fact there exist a diffeomorphism-germ 
$\sigma : (\R^2, 0) \to (\R^2, 0)$ and $\Gamma : (\R^3, 0) \to (\R^3, 0)$ such that 
$\Gamma\circ \widetilde{G} \circ \sigma^{-1} = G$, where $G(t, u) = (u, T_0, T_1)$, the normal form of the swallowtail. 
Moreover it is known that $\sigma$ can be taken to preserve $\{ t = 0\}$, which corresponds to the tangent line 
to the base point of the curve geometrically (\cite{Ishikawa93, Ishikawa95}). 
Then $\sigma$ preserves both ${\mathcal R}_g$ and $t^{i+1}{\mathcal E}_2$ for any $i$, and therefore it preserves 
${\mathcal R}_g^{(i+1)} = {\mathcal R}_g \cap t^{i+1}{\mathcal E}_2$. 

Define $\widetilde{\Gamma} : (\R^4, 0) \to (\R^4, 0)$ by 
$\widetilde{\Gamma}(x_1, x_2, x_3, x_4) = (\Gamma(x_1, x_2, x_3), x_4)$. 
Then $\widetilde{F}$ is diffeomorphic via $(\sigma, \widetilde{\Gamma})$ to 
$H : (\R^2, 0) \to (\R^4, 0)$ defined by $(u, T_0, T_1, h(t, u))$ where $h = (T_3 + \Phi_4(t, u))\circ\sigma^{-1}$. 
Then we have $h \in {\mathcal R}_g^{(4)}$. 
Now, by Lemma \ref{algebraic-lemma}, ${\mathcal R}_g^{(4)}$ is generated by 
$T_3, T_4, T_5$ over $g^*$. Therefore $h$ is expressed as 
$$
h(t, u) = \varphi(u, T_0)T_3 + \psi(u, T_0)T_4 + \rho(u, T_0)T_5, 
$$
for some function-germs $\varphi, \psi, \rho : (\R^2, 0) \to \R$. By Lemma \ref{algebraic-lemma-2}, we have 
$h(t, u) = \Lambda(u, T_0, T_1) + \psi(u, T_0)T_4$, where 
$\Lambda = \varphi(u, T_0)(\frac{1}{2} T_0^2 - uT_1) + \rho(u, T_0)(- \frac{1}{3}uT_0^2 + \frac{2}{3}u^2T_0 + \frac{2}{3}T_1^2)$. 
By the diffeomorphism $(x_1, x_2, x_3, x_4) \mapsto (x_1, x_2, x_3, x_4 - \Lambda(x_1, x_2, x_3))$, 
we have that $\Tan(f)$ is diffeomorphic to the form 
$(t, u) \mapsto ( u, \ T_0, \ T_1, \ \psi(u, T_0)T_4)$ as required. 
If $\psi(0, 0) \not= 0$, then by the diffeomorphism 
$(x_1, x_2, x_3, x_4) \mapsto (x_1, x_2, x_3, \psi(x_1, x_2)^{-1}x_4)$, we have that $\Tan(f)$ is diffeomorphic 
to ${\mbox{\rm USW}}_{2,4}$. 
\QED

\bee
{\rm 
Let $f_\lambda : (\R, 0) \to (\R^4, 0)$ be the family of curves defined by 
$$
f_\lambda(t) = (x_1(t), x_2(t), x_3(t), x_4(t)) = (t^2, \ t^3, \ t^4, \ t^6 + \lambda t^7). 
$$
A tangent frame is given by $\tau_\lambda(t) = (1, \frac{2}{3}t, 2t^2, 3t^4 + \frac{7}{2}\lambda t^5).$ 
Then, by setting $u = -3(t^2+s)$, we see that 
the tangent map $F_\lambda(t, s) = \Tan(f_\lambda)(t, s) = f(t) + s\tau_\lambda(t)$ is diffeomorphic to 
$\widetilde{F}_\lambda : (\R^2, 0) \to (\R^4, 0)$ defined by 
$$
(t, u) \mapsto \left( u, \ t^3+ut, \ \frac{3}{4}t^4 + \frac{1}{2}ut^2, \ \frac{1}{2}t^6+\frac{1}{4}ut^4+\frac{35}{24}\lambda(\frac{3}{7}t^7+\frac{1}{5}ut^5)\right). 
$$
By eliminating $\frac{1}{2}t^6+\frac{1}{4}ut^4$ by the relation 
$\frac{1}{2}t^6+\frac{1}{4}ut^4 = 
\frac{1}{2}(t^3 + ut)^2 - u(\frac{3}{4}t^4 + \frac{1}{2}ut^2)$,  
we have that the germ $F_\lambda$ turns to be diffeomorphic to 
$$
(t, u) \mapsto \left( u, \ t^3+ut, \ \frac{3}{4}t^4 + \frac{1}{2}ut^2, \ \lambda(\frac{3}{7}t^7+\frac{1}{5}ut^5)\right). 
$$
If $\lambda = 0$, then we see that $F_0$ is an embedded swallowtail and is not injective. 
If $\lambda \not= 0$, then we have 
that $F_\lambda$ is diffeomorphic to ${\mbox{\rm generic USW}}_{2,4}$, 
which is injective and therefore $F_\lambda = \Tan(f_\lambda), \lambda \not= 0$ is not diffeomorphic to $F_0 = \Tan(f_0)$. 
}
\ene

\ber
{\rm 
Note that the normal form of ${\mbox{\rm generic USW}}_{2,4}$ is diffeomorphic to the tangent surface to 
the curve of type $(2, 3, 4, 7), t \mapsto (t^2, t^3, t^4, t^7)$. 
The function $\psi$ of Theorem \ref{unfurled-swallowtail-classification}, 
in a certain sense, controls the partial opening of the swallowtail in $\R^4$. 
The geometric meaning of $\psi$ and the exact classification of singularities of the tangent surfaces to 
curves of type $(2, 3, 4, 6)$ for non-generic cases are still open, 
as far as the author knows. The classification for the cases $(2, 3, 4, 2m+1), m \geq 3$ seems to be involved to  
the problem. 
}
\enr

\

\noindent 
ISHIKAWA Goo, \\ 
Department of Mathematics, Hokkaido University, \\
Sapporo 060-0810, JAPAN. \\ 
e-mail : ishikawa@math.sci.hokudai.ac.jp 

\end{document}